\documentclass[11 pt, psamsfonts]{amsart}
\usepackage{amssymb}
\usepackage{array,calc,url,tabularx, epic,tikz}
\usetikzlibrary{arrows,positioning,decorations.pathmorphing,
  decorations.markings}

\def\ignore#1{\relax}

\def\g{\mathfrak g}

\def\sp{{\mathfrak sp}}

\def\R{{\mathbb R}}
\def\Z{{\mathbb Z}}
\def\nat{{\mathbb N}}
\def\Q{{\mathbb Q}}
\def\C{{\mathbb C}}

\def\de{\delta}
\def\la{\lambda}
\def\La{\Lambda}

\def\N{\mathbb N}

\def\Cl{\mathcal C}

\def\P{\mathcal P}
\def\Pnla{{\mathcal P}(n,\la)}

\def\U{{\bf U}}

\def\v{\vskip 2.5mm}

\def\ignore#1{\relax}

\def\om{\omega}
\def\eps{\epsilon}
\def\1{{\bf 1}}
\def\lan{\langle}
\def\ra{\rangle}

\def\ep{\epsilon}

\def\End{{\rm End}}
\def\Hom{{\rm Hom}}

\def\eps{\varepsilon}

\calclayout

\renewenvironment{proof}[1][\proofname]{\par
  \normalfont
  \topsep6\p@\@plus6\p@ \trivlist
  \item[\hskip\labelsep\bfseries
    #1\@addpunct{.}]\ignorespaces
}{%
  \qed\endtrivlist
}

\def\th@plain{%
  \let\thmhead\thmhead@plain \let\swappedhead\swappedhead@plain
  \thm@preskip.5\baselineskip\@plus.2\baselineskip
                                    \@minus.2\baselineskip
  \thm@postskip\thm@preskip
  \itshape
\renewcommand{\labelenumi}{{(\alph{enumi})\quad}}
                        \renewcommand{\labelenumii}{{(\roman{enumii})\ }}
}
\def\th@definition{%
  \let\thmhead\thmhead@plain \let\swappedhead\swappedhead@plain
  \thm@preskip.5\baselineskip\@plus.2\baselineskip
                                    \@minus.2\baselineskip
  \thm@postskip\thm@preskip
  \upshape
}
\def\th@remark{%
  \thm@headfont{\itshape}
  \let\thmhead\thmhead@plain \let\swappedhead\swappedhead@plain
  \thm@preskip.5\baselineskip\@plus.2\baselineskip
                                    \@minus.2\baselineskip
  \thm@postskip\thm@preskip
  \upshape
}

{\theoremstyle{plain}
\newtheorem{theorem}{Theorem}[section]
}

{\theoremstyle{plain}
\newtheorem{proposition}[theorem]{Proposition}
}

{\theoremstyle{plain}
\newtheorem{corollary}[theorem]{Corollary}
}

{\theoremstyle{plain}
\newtheorem{lemma}[theorem]{Lemma}
}

{\theoremstyle{plain}

}

{\theoremstyle{definition}

}

{\theoremstyle{definition}
\newtheorem{example}[theorem]{Example}
}

{\theoremstyle{remark}
\newtheorem{remark}[theorem]{Remark}
}

{\theoremstyle{remark}

}

\numberwithin{equation}{section}

\renewcommand{\labelenumi}{{ \theenumi.}}
\renewcommand{\labelenumii}{{(\alph{enumii})}}

\def\v{\vskip 2.5mm}

\def\la{\lambda}
\def\al{\alpha}

\def\choose #1 #2{\begin{pmatrix}#1\\#2\end{pmatrix}}

\begin{document}

\title[Affine $G_2$ Centralizer Algebras ]
{Affine $G_2$ Centralizer Algebras}

\author{Lilit Martirosyan and Hans Wenzl}

\address{Department of Mathematics\\ University of California\\ San Diego,
California}

\email{hwenzl@ucsd.edu}

\begin{abstract}
We show that $\End_\U(V_\la\otimes V^{\otimes n})$ is generated by the
affine braid group $AB_n$ where $\U=U_q\g(G_2)$, $V$ is its 7-dimensional
irreducible representation and $V_\la$ is an arbitrary irreducible
representation.
\end{abstract}
\maketitle

It is well-known that the famous Schur-Weyl duality between the
general linear group $Gl(N)$ and the symmetric group $S_n$ extends
to a duality between quantum groups of classical Lie types
and the braid group $B_n$, acting on
tensor powers of the vector representations (see \cite{Ji}, \cite{Re}).
More recently, this has also been shown
for the 7-dimensional representation $V$ of Lie type $G_2$. More
precisely, if $\U=U_q\g(G_2)$ is the Drinfeld-Jimbo quantum group,
then $\End_\U(V^{\otimes n})$ is generated by the image of
the braid group $B_n$ acting on $V^{\otimes n}$ via $R$-matrices;
this was shown by Lehrer and Zhang \cite{LZ} and Morrison \cite{Ms},
using earlier results  by Kuperberg
\cite{Ku1}, \cite{Ku2} and Schwarz \cite{schw}.
If $V$ is the smallest nontrivial representation of Lie type $E_N$,
$N=6,7$, a similar results holds as well; one only needs to add
one additional generator in the $(N-1)^{st}$ tensor power
to the $R$-matrices, see \cite{Wexc}.

In the current paper we consider the decomposition of tensor powers
$V_\la\otimes V^{\otimes n}$, where $V_\la$ is an arbitrary irreducible
representation of $U_q\g(G_2)$. It is well-known that the $R$-matrix
formalism now allows a representation of the affine braid group
$AB_n$ into $\End_\U(V_\la\otimes V^{\otimes n})$; here $AB_n$
is the braid group for the Coxeter graph $B_n$.
The main result
of this paper is to show that this map is surjective.
In particular, we retrieve the aforementioned results in the special case
where $V_\la$ is the trivial representation. Besides its intrinsic interest, these results should also be useful
for  categorifying the Lie algebra of type $G_2$. We learned about this from
Catharina Stroppel, who suggested working on these questions.

Our approach is quite different from the ones in the previously mentioned
papers. It is based on the well-known (quantum) Jucys-Murphy
approach and its generalizations in \cite{lr} and \cite{Wexc}.
We expect that it might be useful also for studying at least
parts of tensor powers of adjoint representations of exceptional Lie types.

In more detail, we review necessary tools from the study of quantum groups
in the first section. We then study the case involving the 7-dimensional
representation of $G_2$ in the second section. In particular, we prove
the result mentioned in the abstract. In the third section, we use our
approach to calculate certain structure coefficients in some of the
relations in our algebra. We then close with a discussion of
using the approach in this paper for other Lie types.

$Acknowledgements:$ L. M. would like to thank the Max Planck Institute for Mathematics, Bonn and H.W. would like to thank the Hausdorff Institute, Bonn for hospitality and support while this paper was finished. We would also like to thank Catharina Stroppel for her suggestion to work on this problem.

\section{Quantum groups}

\subsection{Littelmann paths}\label{Litpath}
Let $\g$ be a semisimple Lie algebra, and let $\U=U_q\g$ be the corresponding Drinfeld-Jimbo quantum group.
The reader not familiar with $U_q\g$ should be able to read
this section by just replacing $U_q\g$ by $\g$.
Let $\la$ be a dominant integral highest
weight of $\g$, and let $V_\la$ be the corresponding simple representation
of $U_q\g$.
We also assume that $V$ is a simple representation of $\U$
all of whose weights have multiplicity 1. This allows us to give a
fairly simple description of $\End_\U(V_\la\otimes V^{\otimes n})$ via paths.
This can be viewed as a slight variation of a special case of
Littelmann paths (see \cite{Li});
the latter formalism could be used to give similar
descriptions for arbitrary $V$. The simpler version here has been
known much longer, see e.g. \cite{St} and references there.

For given dominant integral weight $\la$ and $n\in\nat$, we define
the set $\Pnla$ to be the set of all paths
$$t: \la=\la^{(0)}=t(0)\to\la^{(1)}=t(1)\to\ ...\ \to\la^{(n-1)}=t(n-1)\to\la^{(n)}=t(n),$$
where the $\la^{(i)}$'s are dominant integral weights such that
$$V_{\la^{(i+1)}}\subset V_{\la^{(i)}}\otimes V.$$
Note that by our assumptions on $V$, the module $V_{\la^{(i+1)}}$
appears at most with multiplicity one in $V_{\la^{(i)}}\otimes V$.
If $V$ contains the zero weight, two consecutive weights $\la^{(i)}$
and $\la^{(i+1)}$ may coincide. The following theorem is well-known,
and it can be easily proved by induction on $n$.

\begin{theorem}\label{ldecomposition} We have a direct sum
decomposition of $\U$-modules given by
$$V_\la\otimes V^{\otimes n}=\bigoplus_\mu m(\mu,n)V_\mu,$$
where the multiplicity $m(\mu,n)$ is given by the number of paths
in $\Pnla$ which end in $\mu$. In particular, we have
$$ \Cl(n,\la)=\End_\U(V_\la\otimes V^{\otimes n})\ \cong\ \bigoplus_\mu M_{m(\mu, n)},$$
where $M_k$ are the $k\times k$ matrices.
\end{theorem}

\begin{remark}\label{ldecomremark} Let $t\in\Pnla$. Then we denote
by $t'$ the path in $\P(\la, n-1)$ obtained by removing $\la^{(n)}$
\end{remark}

\begin{corollary}\label{pathc}
There exists an assignment $t\in \Pnla\mapsto p_t\in \Cl(n,\la)=
\End_\U(V_\la\otimes V^{\otimes n})$
such that $p_tV^{\otimes n}$ is an irreducible $\g$-module with
highest weight $t(n)$, and such that $p_tp_s=\delta_{ts}p_t$.
The idempotents $p_t$ are uniquely defined by the properties above
and the following one: If $s\in \P(n-1,\la)$, we have
$$p_s=\sum_{t,\ t'=s} p_t.$$
\end{corollary}

One checks easily that $z_\mu^{(n)}=\sum_{t\in \Pnla(\mu)}p_t$ is
a central idempotent in $\Cl(n,\la)=\End_\U(V_\la\otimes V^{\otimes n})$.
 Consider the subalgebra
$\Cl(n-1,\la)\otimes 1\subset \Cl(n,\la)$;
if no confusion arises we will usually
denote the latter algebra only by $\Cl(n-1,\la)$. Let $W^{(n)}_\nu$ be
a simple $\Cl(n,\la)$-module labeled by the dominant weight $\nu$.
Then we have the following isomorphism of $\Cl(n-1,\la)$-modules:
\begin{equation}\label{restriction}
W^{(n)}_\nu \cong \oplus_\mu  W^{(n-1)}_\mu,
\end{equation}
where $\mu$ runs through all highest weights in $V^{\otimes (n-1)}$
such that $V_\nu\subset V_\mu\otimes V$. Keeping with the notation
of Section \ref{Litpath}, we may also refer to this situation as
$\nu$ being connected with $\mu$ by a path of length 1.

By definition, we can define a basis $(v_t)_{t\in\P(n,\la)_\nu}$
for the simple $\Cl(n,\la)_\nu$-module $W^{(n)}_\nu$
labeled by all paths of length $n$ in $\Pnla$ which end in $\nu$.
here the vector $v_t$ spans the image
of $p_t$ for each $t\in \Pnla(\nu)$ and it is uniquely determined
up to scalar multiples.
Let $\delta, \nu$ be dominant weights for which
$V_\delta\subset V^{\otimes n-k}$ and $V_\nu\subset V^{\otimes n}$,
and let $\P_k(\delta,\nu)$ be the set of all paths of length $k$
from $\delta$ to $\mu$, with paths as defined
in Section \ref{Litpath}. Let $W(\delta,\nu)$ be the vector space
spanned by these paths. Then we obtain a representation of
$ \Cl(k,\la)=\End_\U(V_\la\otimes V^{\otimes k})$ on $W(\delta,\nu)$ by
\begin{equation}\label{skew}
a\in \End(V^{\otimes k}) \mapsto (p_t\otimes a)z_\nu^{(n)};
\end{equation}
here we used the obvious bijection between elements $s\in \P_k(\delta,\nu)$
and paths $\tilde s\in \P_n(\nu)$
for which $\tilde s_{|[0,n-k]}=t$.

\subsection{Generating $\End_\U(V_\la\otimes V^{\otimes n})$}\label{generatingEnd}
In order to find a generating set for $\End_\U(V_\la\otimes V^{\otimes n})$
we will use a simple lemma as follows. Fix a dominant weight $\nu$
such that $V_\nu\subset V_\la\otimes V^{\otimes n}$. Let $\mu_1,\mu_2$ be
dominant weights such that $W_{\mu_i}^{(n-1)}\subset W_\nu^{(n)}$
for $i=1,2$, with notations as in Eq \ref{restriction}.
We say that $\mu_1$ is equivalent to $\mu_2$ if there exists a dominant
integral weight $\delta$ such that $V_{\mu_i}\subset V_\delta\otimes V$
for $i=1,2$ and $V_\delta\subset V^{\otimes n-2}$, and an element
$a\in\End_\U(V^{\otimes 2})$ such that all the matrix coefficients of
$a$, acting on the path space labeled by the elements of $\P_2(\delta,\nu)$
are nonzero.

\begin{lemma}\label{generatinglemma} Assume that for given $W_\nu^{(n)}$
all the dominant weights
$\mu$ in Eq \ref{restriction} are equivalent in the sense just defined.
Then the algebra $\Cl_n'$ generated by $\Cl_{n-1}$ and $1_{V_\la\otimes V^{\otimes n-2}}\otimes \End_\U(V^{\otimes 2})$ acts irreducibly on
$W_\nu^{(n)}$. In particular, if this holds for all $\nu$ for which
$V_\nu\subset V_\la\otimes V^{\otimes n}$, then the algebra $\Cl_n'$ coincides
with $\Cl_n$.
\end{lemma}

$Proof.$ Assume $W$  is a nonzero $\Cl_n'$ submodule of $W_\nu^{(n)}$.
Then it suffices to show that it must coincide with the right hand side
of Eq \ref{restriction} if viewed as a $\Cl_{n-1}$ module. Assume to the
contrary that there exists a $\mu_1$ such that $W^{(n-1)}_{\mu_1}\not\subset W$. As $\mu_1$ is equivalent to any $\Cl_{n-1}'$ submodule of $W$, by assumption, we can find a $\mu_2$ such that $W^{(n-1)}_{\mu_2}\subset W$,
a dominant weight $\delta$ and an $a\in \End_\U(V^{\otimes 2})$
which satisfy
the conditions of equivalence, as stated at the beginning of this section.
But then, in particular, $1_{V_\la\otimes V^{\otimes n-2}}\otimes a$
does not map $W_{\mu_2}^{(n-1)}$ into $W$; just apply it to a basis vector
$v_t$ for which $t(n-1)=\mu_1$ and $t(n-2)=\delta$. Hence $W$ is not
a $\Cl_n'$-module.
This proves the first claim of the statement. This also implies the
second statement, as we have a faithful representation of $\Cl_n$ on
the direct sum of all $W_\nu^{(n)}$ labeled by the $\nu$s for which
$V_\nu\subset V_\la\otimes V^{\otimes n}$.

\subsection{Quantum groups}\label{quantum}
 We assume as ground ring the field
$\Q(q)$ of rational functions in the variable $q$; most of the results
hold in greater generality (e.g. if $q$ is a complex number not equal to
a root of unity or 0). It is well-known
that in our setting the category $Rep(\U)$
of integrable representations of $\U$ is semisimple,
and it has the same Grothendieck semiring
as the original Lie algebra. Moreover, $Rep(\U)$ is a braided tensor
category. This implies that for $\U$-modules $V$, $W$,
there are natural braiding isomorphisms
$ R_{VW}:V\otimes W\rightarrow W\otimes V$ which satisfy
\begin{equation}\label{braiding}
R_{U, V\otimes W}=(1_V\otimes R_{UW})(R_{UV}\otimes 1_W),
\end{equation}
where $U,V,W$ are $\U$-modules; a similar formula holds for
$R_{U\otimes V,W}$.
Let $B_n$ be Artin's braid groups, given by generators $\sigma_i,
\ 1\leq i\leq n-1$ and relations $\sigma_i\sigma_{i+1}\sigma_i=
\sigma_{i+1}\sigma_i\sigma_{i+1}$ as well as $\sigma_i\sigma_j=\sigma_j
\sigma_i$ for $|i-j|\geq 2$.
Moreover, let $AB_n$ be the {\it affine braid group}, where we add
to the generators of $B_n$ the additional generator $\tau$, and the
additional relations
$$\sigma_1\tau\sigma_1\tau=\tau\sigma_1\tau\sigma_1\quad{\rm and}\quad
\sigma_i\tau=\tau\sigma_i,\quad {\rm for\ i>1}.$$
Equivalently, $AB_n$ is the braid group of Dynkin type $B_n$;
we shall not use this notation again, so $B_n$ will denote the braid group
in the rest of this paper.
We obtain, for any  $\U$-module $V$,
a representation of Artin's braid group $B_n$ in $\End(V^{\otimes n})$
by the map
$$\sigma_i\mapsto R_i=1_{i-1}\otimes R_{VV}\otimes 1_{n-1-i}
\in \Cl_n=\End_\U(V^{\otimes n}),$$
where $1_j$ is the identity map on $V^{\otimes j}$.
Similarly, if $V_\la$ is a simple highest weight module of $\U$,
we can extend the just defined representation of $B_n$ to a representation
of $AB_n$ in $\End_\U(V_\la\otimes V^{\otimes n})$ by
$$\tau\ \mapsto R_{V,V_\la}R_{v_\la,V}\otimes 1_{n-1}\ \in\ \End_\U(V_\la\otimes V^{\otimes n}),$$
and where the action of the $\sigma_i's$ is given by the previously
defined representation with the obvious embedding of
$\End_\U(V^{\otimes n})$
into $\End_\U(V_\la\otimes V^{\otimes n})$.

As we also have a faithful representation of $\Cl(n,\la)$
with respect to the basis $(v_t)$, with $t\in\Pnla$, we also
obtain a representation of $AB_n$ on this path basis.
In particular, we obtain matrices $A_i$ such that
\begin{equation}\label{pathrep}
\sigma_i\mapsto A_i: t\to \sum_s a_{st}^{(i)} s.
\end{equation}
If $t\in \Pnla(\mu)$, then so are the paths $s$ in the equation above.
Moreover, it follows from Eq. \ref{skew} that
the paths $s$ differ from $t$ only in the interval $[i-1,i+1]$.
Because of this we shall often only consider
spaces $W_i(\la,\nu)$ with a basis consisting of paths of length 2
from $\la=t(i-1)$ to $\nu=t(i+1)$. Equation \ref{pathrep} induces an
obvious action of $A_i$ on $W_i(\la,\nu)$. We will call the corresponding
matrix block $A_i(\la,\nu)$. If there is no danger of confusion,
we will often suppress the index $i$ in $W_i(\la,\nu)$ and  $A_i(\la,\nu)$.
We shall also need the following theorem,
due to Drinfeld \cite{Dr}.
\begin{proposition}
\label{prop:eigenvalue}
Let $V_\la, V_\mu, V_\La=V$ be
simple $\U$-modules with highest weights
$\la, \mu, \La$ respectively, and such that $V_\mu$ is a submodule of
$V_\la\otimes V_\La$. Then
$$(R_{V_{\la}V_{\La}}R_{V_{\La}V_{\la}})_{|V_\mu}=
q^{c_\mu- c_\la-c_\La}1_{V_\mu},$$
where for any weight $\gamma$ the quantity $c_\gamma$ is given by
$\lan\gamma +2\rho ,\gamma\ra$.
\end{proposition}
\noindent

\subsection{Representations of $AB_2$}\label{AB2section} In the following we consider
representations of the affine braid group $AB_2$ on a finite dimensional vector space $W$. In more detail,
we consider matrices $A$ and $T$ which act on $W$ and which satisfy
 the following conditions:
\vskip .2cm
\renewcommand{\labelenumi}{\alph{enumi})}
\begin{enumerate}
\item They satisfy the braid relation $ATAT=TATA$,
\item The matrix $A$ satisfies the relation $A-A^{-1}=(q-q^{-1})(1-mP)$, where $m=(r-r^{-1} )/(q-q^{-1})$, and where $P$ is a rank 1 eigenprojection of $A$.
\item The central element $TATA$ acts as the identity on $W$.
\item We assume that $T$ is a diagonal matrix with eigenvalues $q^{e(t)}$
    where $t$ runs through a labeling set for a basis of $W$.
\end{enumerate}

\vskip .3cm
\begin{proposition}\label{AB2representations}
(a) The matrix entries of $A$ and $P$ are related by the equation
$$(1-q^{e(t)+e(s)})a_{ts}=(q-q^{-1})\delta_{ts}-(r-r^{-1}+q-q^{-1})p_{ts}.$$

(b) The diagonal entry $d_s=p_{ss}$ is equal to zero only if $e(s)=\pm 1$.
\end{proposition}

$Proof.$ The proofs of these statements are variations of the proofs
of \cite{Wexc} Lemma 4.1 and Lemma 4.3. For (a), just observe that by our
assumption (c) we have $A^{-1}=TAT$. The claim now follows by plugging
this expression for $A^{-1}$ into assumption (b).

To prove part (b), we write the rank 1 idempotent $P$ in the form
$P=vw^T$ for appropriate column vectors $v$ and $w$. Multiplying the equation in assumption (b) by $P$, we deduce that the eigenvalue corresponding to $P$
is equal to $r^{-1}$. It follows that $Av=r^{-1}v$ and $w^TA=r^{-1}w^T$.
If $p_{ss}=0$ for some index $s$, it follows that the $s$-th row or the $s$-th column of $P$ is equal to zero. The same applies to $A$ except
for the diagonal entry $a_{ss}$, by (a). But then $a_{ss}$ must be an
eigenvalue of $A$. If $a_{ss}=q$, we deduce from (a) that $e(s)=-1$.
Similarly, $a_{ss}=-q^{-1}$ implies that $e(s)=1$. Finally, if $a_{ss}=r^{-1}$, then the $s$-th standard basis vector would be a left
eigenvector (or its transpose a right eigenvector) of $A$ which would
be different from $v$ resp. $w$. This would contradict the fact that
the eigenvalue $r^{-1}$ has multiplicity 1.

\begin{corollary}\label{irred} Assume that the eigenvalues $q^{e(t)}$ of $T$ are mutually distinct, and none of them is equal to $q^{\pm 1}$.
Then the representation of $AB_2$ on $W$ is indecomposable.
\end{corollary}

$Proof.$ It suffices to show that $A$ and $T$ generate the full matrix
algebra over $W$. We get the diagonal matrix entries $E_{ss}$ as eigenprojections of $T$, and the off-diagonal matrix entries $E_{st}$ as $\frac{1}{ p_{st}}E_{ss}PE_{tt}$.
As $P$ is a rank 1 idempotent,
its entries $p_{st}$ and $p_{ts}$ are nonzero
because $p_{st}p_{ts}=p_{ss}p_{tt}\neq 0$.

\begin{lemma}\label{AB2rep} ($q$-Murphy-Jucys approach)
Assume now that $AB_2$ acts on a two-dimensional vector space $W$  on
which $A$ satisfies the relation $A-A^{-1}=(q-q^{-1})1$.
Moreover, $T$
has the two distinct eigenvalues $q^{\pm e(t)}$ such that $TATA=1$.
Then $A$ has non-zero off-diagonal matrix entries, except possibly
if $T=1$.
\end{lemma}

$Proof.$ It follows as in the proof of Proposition \ref{AB2representations} that
$$(1-q^{e(t)+e(s)})a_{ts}=(q-q^{-1})\delta_{ts}.$$
By assumption, we have $e(s)+e(t)=0$
for the two basis paths $s$ and $t$, and $e(s)\neq 0$.
It follows that $a_{ss}=(q-q^{-1})/(1-q^{2e(s)})$.
A similar formula also holds for $a_{tt}$. If $P$ now is the
eigenprojection of $A$ for, say, $q$, we have $P=(A+q^{-1}1)/(q+q^{-1})$.
We deduce from this that the diagonal entries of $P$ are nonzero and,
as $P$ is a rank 1 idempotent, so are also its off-diagonal entries.
As $A=qP-(q+q^{-1})1$, the same also holds for the off-diagonal entries of
$A$.

\subsection{}\label{quantumrep}
Let now $\U=U_q\g$ be a Drinfeld-Jimbo quantum group, and let
$V$, $V_\la$ be simple $U_q\g$-modules, with $\la$ being the highest weight
for $V_\la$. We would like to apply the results from the previous section
for certain representations of $AB_2$ appearing in $\End_\U(V_\la\otimes V^{\otimes 2})$. More precisely, let $V_\nu\subset V_\la\otimes V^{\otimes 2}$ be an irreducible $U_q\g$ module. Observe that
$W=Hom_\U(V_\nu, V_\la\otimes V^{\otimes 2})$ is a vector space whose dimension is equal to the multiplicity of $V_\nu$ in
$V_\la\otimes V^{\otimes 2}$. Using notations of Section \ref{quantum},
we can make $W$ into an $AB_2$-module by mapping
$$A\mapsto \alpha^{-1} (1\otimes R_{V,V}),\hskip 3em
T\mapsto \gamma (R_{V,V_\la}R_{V_\la,V}\otimes 1),$$
acting via concatenation of morphisms on $W$; here $\alpha$ and $\gamma$
are scalars which will be fixed below as follows:

1. We assume that $\nu$ is such that $A$ has at most three distinct
eigenvalues. Moreover, $\alpha$ is chosen such that two of these eigenvalues are $q$ and $-q^{-1}$. The third eigenvalue, denoted by $r^{-1}$
has multiplicity 1 in $W$.

To check this in practice, assume $V^{\otimes 2}=\oplus_\kappa V_\kappa$.
Then, for a fixed $\kappa$, the eigenprojection of $A$ corresponding to $\kappa$ has multiplicity 1 on $W$ if and only if $V_\nu$ appears with multiplicity 1
in $V_\la\otimes V_\kappa\subset V_\la\otimes V^{\otimes 2}$.

2. We assume $\gamma$ is chosen such that $TATA=1$.

\begin{lemma}\label{normalizing} Let $V_\la\otimes V=\bigoplus_t V_{\mu_t}$,
and fix a dominant weight $\nu$ such that $V_\nu\subset V_\la\otimes V^{\otimes 2}$.
Then there exist scalars $\alpha$ and $\gamma$ such that the assumptions above are satisfied, and such that the eigenvalues of $T$ are given by $q^{e(t)}$, where
$$2e(t)= 2c_{\mu_t}-c_\la-c_\nu+2e_\alpha;$$
here $e_\alpha$ is given by $\alpha=q^{e_\alpha}$, and $c_\kappa$ is given by
$$c_\kappa=(\kappa+2\rho,\kappa)$$
for any weight $\kappa$, with $\rho$ being half the sum of the positive roots
of $\g$.
\end{lemma}

$Proof.$ This is a straightforward consequence of properties of $R$-matrices. Let $A'=1\otimes R_{V,V}$ and $T'=R_{V,V_\la}R_{V_\la,V}\otimes 1$,
acting on $V_\la\otimes V^{\otimes 2}$. By Prop \ref{prop:eigenvalue},
$T'$ acts on the vector $v_t$ corresponding to the path
$t: \delta\to\mu_t\to\nu$ by the scalar $q^{c(t)-c_V-c_\delta}$,
where $c(t)$ and $c_V$ are the values of the Casimir for $V_{\mu_t}$ and $V$. By the braiding axiom \ref{braiding} the element
$A'T'A' = R_{V^{\otimes 2},V_\la}R_{V_\la,V^{\otimes 2}}$
acts on $V_\nu\subset V_{\mu_t}\otimes V\subset V_\la\otimes V^{\otimes 2}$ by the scalar $q^{c_\nu-c(t)-c_V}$.
It follows that $T'A'T'A'$ acts by the scalar $q^{c_\nu-c_\la-2c_V}$.
So if $A=q^{-e_\al}A'$ and $T=q^{e_\al+c_V+(c_\de-c_\nu)/2}T'$,
we get the desired identity $TATA=1$. The values of the eigenvalues of
$T$ follow from the formulas in this proof.

\subsection{Rough gradation of tensor products} We assume $V$ to be a
self-dual representation of a quantum group. Using semisimplicity
of our representation category, we decompose the tensor product $V_\la\otimes V^{\otimes n}$ as
$$V_\la\otimes V^{\otimes n}= (V_\la\otimes V^{\otimes n})_{old}\oplus (V_\la\otimes V^{\otimes n})_{recent}\oplus
(V_\la\otimes V^{\otimes n})_{new};$$
here $(V_\la\otimes V^{\otimes n})_{new}$ is the maximum direct sum of simple representations which have not appeared in any lower tensor product,
$(V_\la\otimes V^{\otimes n})_{recent}$ is the maximum direct sum of simple representations which have appeared in the $V_\la\otimes V^{\otimes n-1}$ for the first time,
and $(V_\la\otimes V^{\otimes n})_{old}$
is the maximum direct sum of simple representations which have already appeared in $V_\la\otimes V^{\otimes k}$ for $k<n-1$. We have the following result
(see e.g. \cite{Wexc}, Prop. 4.10; the proof there for the case $\la=0$ works as well in this more general setting):

\begin{proposition}\label{tensorold}
The algebra $\End_\U(V_\la\otimes V^{\otimes n})_{old}$ is generated by
the restrictions to $(V_\la\otimes V^{\otimes n})_{old}$ of $\End_\U(V_\la\otimes V^{\otimes n-1})\otimes 1$ and $1_{V_\la\otimes V^{\otimes n-2}}\otimes Q$, where $Q$ is the projection onto the trivial representation $\1\subset V^{\otimes 2}$.
\end{proposition}

\begin{remark} The proof uses Lemma \ref{generatinglemma} by showing that the matrix coefficients of the projection $Q$ for suitable path bases are
nonzero.
\end{remark}

\section{The example $G_2$}

\subsection{Example} We will be particularly interested in the case
with $\g=\g(G_2)$ and $V$ its simple 7-dimensional representation.
We first recall some basic facts about its roots and weights.

With respect to the orthonormal unit vectors $\eps_1$, $\eps_2$, $\eps_3$ of $\R^3$, the roots of $\g$ can be written $\Phi=\pm \{\eps_1-\eps_2,\eps_2-\eps_3, \eps_1-\eps_3, 2\eps_1-\eps_2-\eps_3, 2\eps_2-\eps_1-\eps_3, 2\eps_3-\eps_1-\eps_2\}$. The base can be chosen $\Pi=\{\al_1=\eps_1-\eps_2, \al_2=-\eps_1+2\eps_2-\eps_3\}$. The Weyl vector is given by $\rho=2\eps_1+\eps_2-3\eps_3$ and the Weyl group is $D_6$. The fundamental dominant weights are given by $\Delta=\{\La_1=\eps_1-\eps_3, \La_2=\eps_1+\eps_2-2\eps_3 \}$.

The following is the Weyl chamber:

\begin{center}
  \begin{tikzpicture}[scale=.4]
    \draw (0,-3) node[anchor=east]  {Weyl Chambre for $\g$};
   \foreach \y in {0,...,7} \foreach \x in {0,...,\y}
    \draw[xshift=\x cm,thick,blue,fill=blue](\x cm,3*\x cm) circle (.2cm);
       \foreach \x in {0,...,5}
    \draw[xshift=\x cm,thick,blue,fill=blue] (\x cm,3*\x cm +6 cm) circle (.2cm);
        \foreach \x in {0,...,3}
    \draw[xshift=\x cm,thick,blue,fill=blue] (\x cm,3*\x cm +12 cm) circle (.2cm);
     \foreach \x in {0,...,1}
    \draw[xshift=\x cm,thick,blue,fill=blue] (\x cm,3*\x cm +18 cm) circle (.2cm);
    \foreach \x in {2,...,7}
    \draw[xshift=\x cm] (\x cm,3*\x cm) node [anchor=north west]{\tiny $\x\La_1$};
      \foreach \x in {2,...,5}
    \draw[xshift=\x cm] (\x cm,3*\x cm + 6 cm) node [anchor=north]{\tiny $\x\La_1+\La_2$};
     \foreach \x in {2,...,3}
    \draw[xshift=\x cm] (\x cm,3*\x cm + 12 cm) node [anchor=north]{\tiny $\x\La_1+2\La_2$};
        \foreach \y in {2,...,3}
    \draw[yshift=\y cm] (0cm, 5*\y cm) node [anchor=north east]{\tiny $\y\La_2$};
         \foreach \y in {2,...,3}
    \draw[yshift=\y cm] (2.5 cm, 5*\y cm+ 3cm) node [anchor=north]{\tiny $\La_1+\y\La_2$};
   	
	  \draw (0 cm, 0 cm) node[anchor=north west]  {\tiny $0$};
       \draw (0 cm, 6 cm) node[anchor=north east]  {\tiny $\La_2$};
       \draw (2.2 cm, 3 cm) node[anchor=north west]  {\tiny $\La_1$};
       \draw (2.2 cm, 9 cm) node[anchor=north]  {\tiny $\La_1+\La_2$};
 	  \draw[dotted,thick] (0,0) -- +(14, 21);
 	  \draw[dotted,thick] (0,0) -- +(0, 21);
 	 	  \draw[dotted,thick, magenta] (2,3) -- +(0, 18);
 	  \draw[dotted,thick, magenta] (2.1,3.1) -- +(12, 18);


  \end{tikzpicture}
\end{center}

We have, for any $\g$-module $M$,
and for any simple $\g$-module $V_\delta$ that
$$mult_{M\otimes V}\ V_\de=\sum_{\mu} mult_M\ V_\mu,$$
where the sum is over all $\mu$ such that $V_\de\subset V_\mu\otimes V$ with $V=V_{\La_1}$. We will use this simple observation for $M=V_\la\otimes V^{\otimes n-1}$ or $M=V^{\otimes n-1}$.

\begin{proposition} We have $V_\de\subset V_\mu\otimes V$ if $\de=\mu+\om$, where $
\om\neq 0$ is a short root (nonzero weight) of $\g$ or if $\de=\mu$ such that $\de=a\La_1+b\La_2$ with $a\geq 1$.
\end{proposition}

\begin{remark}\label{tensorgeometry} The tensor product rules can be
easily visualized as follows: Consider the hexagon centered at $\mu$
and with corners $\mu+\om$, with $\om$ running through
the short roots of $\g$. If this hexagon is contained in the dominant
Weyl chamber $C$, then $V_\mu\otimes V$ decomposes into the direct sum
of irreducibles $\g$-modules whose highst weights are given
by the corners and the center of
the hexagon. If it is not contained in $C$, leave out all
the corners of the hexagon which are not in $C$; moreover, if
$\mu=b\om_2$, also leave out $\mu$ itself.
\end{remark}

Using the proposition, we can draw the Bratteli diagram for $V^{\otimes n}$.

  \begin{tikzpicture}[scale=.4]

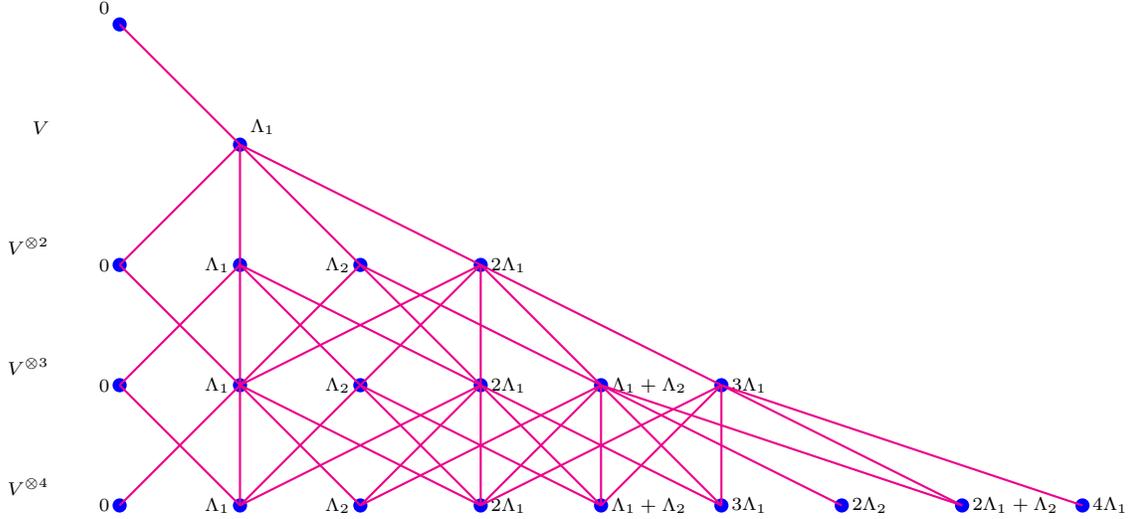
    \draw (0,-16) node[anchor=east]  {Bratteli diagram for $V^{\otimes n}$};
        \draw (-10, 4) node[anchor=south east]  {\tiny $0$};

    \draw [thick,blue,fill=blue](-10 cm,4 cm) circle (.2cm);

     \draw (-6, 0) node[anchor=south west]  {\tiny $\La_1$};
     \draw (-12, 0) node[anchor=south east]  {\tiny $V$};
\foreach \x in {2,...,4}      \draw (-12, -4*\x+4) node[anchor=south east]  {\tiny $V^{\otimes \x}$};

\foreach \x in {1,...,3}     \draw (-10, -4*\x) node[anchor= east]  {\tiny $0$};
\foreach \x in {1,...,3}     \draw (-6, -4*\x) node[anchor= east]  {\tiny $\La_1$};
\foreach \x in {1,...,3}     \draw (-2, -4*\x) node[anchor= east]  {\tiny $\La_2$};
\foreach \x in {1,...,3}     \draw (2, -4*\x) node[anchor= west]  {\tiny $2\La_1$};
\foreach \x in {2,3}     \draw (6, -4*\x) node[anchor= west]  {\tiny $\La_1+\La_2$};
\foreach \x in {2,3}     \draw (10, -4*\x) node[ anchor= west]  {\tiny $3\La_1$};
\draw (14, -12) node[anchor= west]  {\tiny $2\La_2$};
 \draw (18, -12) node[anchor= west]  {\tiny $2\La_1+\La_2$};
  \draw (22, -12) node[anchor= west]  {\tiny $4\La_1$};

    \draw [thick,blue,fill=blue](-6 cm,0 cm) circle (.2cm);
 \foreach \x in {0,...,3}
   \draw[xshift=\x cm,thick,blue,fill=blue] (3*\x cm-10 cm,-4 cm) circle (.2cm);

 \foreach \x in {0,...,5}
   \draw[xshift=\x cm,thick,blue,fill=blue] (3*\x cm-10 cm,-8 cm) circle (.2cm);
\foreach \x in {0,...,8}
\draw[xshift=\x cm,thick,blue,fill=blue, outer sep=0.3pt,
    inner sep=0.5pt] (3*\x cm-10 cm,-12 cm) circle (.2cm);

      \draw[thick, magenta] (-10,4) -- +(4, -4);

    	  \draw[thick, magenta] (-6,0) -- +(8, -8);
    	  \draw[thick, magenta] (-6,0) -- +(16, -8);
    	  \draw[thick, magenta] (-6,0) -- +(0, -12);
	  \draw[thick, magenta] (-10,-4) -- +(4, 4);
	  \draw[thick, magenta] (-6,-4) -- +(8, -4);
	  \draw[thick, magenta] (-6,-4) -- +(4, -4);

	  \draw[thick, magenta] (-6,-8) -- +(8, -4);

    	  \draw[thick, magenta] (-10,-4) -- +(8, -8);
    	  \draw[thick, magenta] (-10,-8) -- +(4, -4);
    	  \draw[thick, magenta] (-10,-8) -- +(4, 4);
    	  \draw[thick, magenta] (-10,-12) -- +(8, 8);
    	  \draw[thick, magenta] (-2,-8) -- +(4, -4);
    	  \draw[thick, magenta] (-2,-8) -- +(8, -4);
    	  \draw[thick, magenta] (-2,-4) -- +(8, -4);
    	  \draw[thick, magenta] (-2,-8) -- +(-4, -4);
    	  \draw[thick, magenta] (2,-4) -- +(0, -8);
    	  \draw[thick, magenta] (2,-4) -- +(-8, -4);
    	  \draw[thick, magenta] (2,-4) -- +(4, -4);
    	  \draw[thick, magenta] (2,-8) -- +(-8, -4);
    	  \draw[thick, magenta] (2,-8) -- +(8, -4);
    	  \draw[thick, magenta] (2,-8) -- +(4, -4);
    	  \draw[thick, magenta] (2,-8) -- +(-4, -4);
    	  \draw[thick, magenta] (2,-4) -- +(-4, -4);
    	  \draw[thick, magenta] (6,-8) -- +(-4, -4);
    	  \draw[thick, magenta] (6,-8) -- +(0, -4);
    	  \draw[thick, magenta] (6,-8) -- +(-8, -4);
    	  \draw[thick, magenta] (6,-8) -- +(4, -4);
    	  \draw[thick, magenta] (6,-8) -- +(8, -4);
    	  \draw[thick, magenta] (6,-8) -- +(12, -4);
    	  \draw[thick, magenta] (10,-8) -- +(0, -4);
    	  \draw[thick, magenta] (10,-8) -- +(12, -4);
    	  \draw[thick, magenta] (10,-8) -- +(-4, -4);
    	  \draw[thick, magenta] (10,-8) -- +(8, -4);
    	  \draw[thick, magenta] (10,-8) -- +(-8, -4);

     \end{tikzpicture}

\begin{example}\label{G2parameters}  We let $V$ be the 7-dimensional
irreducible representation, $V=V_{\La_1}$ with $\La_1$ being
the first fundamental weight. We then get
$$V^{\otimes 2}\cong (\1\oplus V_{2\La_1})\oplus (V\oplus V_{\La_2}).$$
Here the first two summands span the symmetrization, and the third and fourth summand span the antisymmetrization of $V^{\otimes 2}$.
We normalize the invariant product on the weight lattice such that
$(\La_1,\La_1)=2$ and $(\La_2,\La_2)=6$. With these conventions we get
the values $c_\nu=0, 28, 12, 24$ for $\nu=0,2\La_1,\La_1,\La_2$.
Hence the eigenvalues of $R_{V,V}$ are given by $q^{-12}, q^{2},-q^{-6}$
and $-1$. So if $\alpha=q$, $A=\alpha^{-1}R_{V,V}$ has the desired eigenvalues $q$ and $-q^{-1}$ for the representations $V_{2\La_1}$ and $V_{\La_2}$. As we shall see in a moment, it will be convenient to
associate with $P$ the eigenprojection of $A$ projecting onto
$V=V_{\La_1}\subset V^{\otimes 2}$, which corresponds to the eigenvalue $-q^{-7}$. Indeed, let
$W=Hom_\U(V_\nu, V_\la\otimes V^{\otimes 2})$ be as in Section \ref{quantumrep} for $\U=U_q\g(G_2)$. Then $P$ is the projection onto
the subspace $Hom_\U(V_\nu, V_\la\otimes V_{\La_1})$ of $W$,
given by the embedding $V\subset V^{\otimes 2}$.
As all weights of $V$ have multiplicity 1,
the projection $P$ has at most rank 1 in its action on $W$.
Moreover, if $\nu\neq\la$,
$$V_\nu\not\subset V_\la\otimes \1\subset V_\la\otimes V^{\otimes 2}.$$
This shows that the second condition in Section \ref{AB2section} is satisfied.
\end{example}

\subsection{New and recent modules}\label{newrecentsection}
In the following, we assume
$V$ to be as in the last section. Moreover, let $n(\nu)$ be the smallest
integer $n$ for which $V_\nu\subset V_\la\otimes V^{\otimes n}$, for any dominant integral weight $\la$. The following example will be important
for the general case.

\begin{example}\label{nnu} The dominant integral weights for $G_2$
are given by $\nu=(\nu_1,\nu_2,\nu_3)$ with the $\nu_i$ being integers
satisfying $\nu_1\geq\nu_2\geq 0\geq\nu_3$ and $\nu_1+\nu_2+\nu_3=0$.
Then $n(\nu)=\nu_1+\nu_2$, if $V_\la$ is the trivial representation.
Indeed, the sum of the first two coordinates of a weight of $V$ is
at most 1; hence the sum of the first two coordinates of a weight
in $V^{\otimes n}$ is at most $n$. On the other hand, it is easy to
prove by induction on $n$ and the tensor product rule, see Remark \ref{tensorgeometry} and the proposition just before it,
that if
$\nu_1+\nu_2=n$ for an integral dominant weight $\nu$, then $V_\nu$
appears in $V^{\otimes n}$.
\end{example}

\begin{lemma}\label{n(nu)} Let $n_0(\nu)$ be the number as in
Example \ref{nnu}, for $V_\lambda$ being the trivial representation.
Then the number $n(\nu)$ for the general case is given by
$n(\nu)=n_0(w(\nu-\lambda))$, where $w$ is the element of the Weyl group
which maps $\nu-\lambda$ into the dominant Weyl chamber.
\end{lemma}

$Proof.$ Let $\om_1=(1,0,-1)$ and $\om_2=(0,1,-1)$. Then $\nu=\nu_1\om_1+\nu_2\om_2$. The result of the Example \ref{nnu}
can now be phrased  as follows: There exists a path $t$ from 0 to $\nu$
such that $t(i+1)-t(i)\in\{ \om_1,\om_2\}$ and $V_{t(i+1)}\subset
V_{t(i)}\otimes V$. For general $\lambda$, the same approach works
if $\nu-\lambda$ is in the dominant Weyl chamber by just shifting the
path from 0 to $\nu-\lambda$ by $\lambda$. If $\nu-\lambda$ is not dominant,
we first build the path from $\lambda$ to $\lambda+w(\nu-\lambda)$,
with $w$ as in the statement.
We then get the desired path from $\lambda$ to $\nu$ by replacing
the line segments $\om_1$ and $\om_2$ by $w^{-1}(\om_1)$ and
$w^{-1}(\om_2)$, i.e. we reflect the path by $w$, with the fixed point of
the reflection group shifted to $\la$. This shows that $n(\nu)\leq
n(\la+w(\nu-\la))$. On the other hand, if there was a shorter path from
$\la$ to $\nu$, we could reflect it back to a shorter path from $\la$ to $w(\nu-\la)$, contradicting the result already established in this case.

\subsection{Surjectivity} We can now prove one of the main results
of this paper. We will first need the following technical lemma.

\begin{lemma}\label{checke(t)}
The quantity $e(t)$ as in Lemma \ref{normalizing} is not equal to $\pm 1$
if $\delta\neq \nu$ and if there is more than one path from $\delta$ to $\nu$.
\end{lemma}

$Proof.$ Let $V_\nu\subset V_\delta\otimes V^{\otimes 2}$ be a simple subrepresentation,
and let $t$ be the path
$$t: \delta \to \mu_t\to \nu, \quad {\rm with\ } \beta_t=\mu_t-\la,\ \tilde\beta_t=\nu-\mu_t.$$
By Lemma \ref{normalizing} and Example \ref{G2parameters} we have
$$e(t)=c_{\mu_t}-\frac{1}{2}(c_\delta+c_\nu)+1.$$
It follows from a straightforward calculation that
this can also be expressed as
$$e(t)=(\la+\rho, \beta_t-\tilde\beta_t)-(\beta_t,\tilde\beta_t)+\frac{1}{ 2}((\beta_t,\beta_t)-(\tilde\beta_t,\tilde\beta_t))+1.$$
In analyzing this, we will use the following observations:

(A) If $\tilde\beta_t\neq \pm\beta_t$ and both $\beta_t$ and $\tilde\beta_t$ are not equal to 0, then $\beta_t-\tilde\beta_t$ is a (short or long) root.
Depending on whether this root is positive or negative, both
$(\rho,\beta_t-\tilde\beta_t)$ and $(\delta+\rho,\beta_t-\tilde\beta_t)$
have the same sign (including 0), for any dominant weight $\delta$.
If $\delta\neq \nu$, we are left with the following cases:

{\it Case 1:} $\beta_t\neq 0\neq \tilde\beta_t$ and $\beta_t+\tilde\beta_t$ is a short root. In this case $(\beta_t,\tilde\beta_t)=-1$,
$\beta_t-\tilde\beta_t$ is a long root and the formula for
$e(t)$ simplifies to
$$e(t)=(\de+\rho,\beta_t-\tilde\beta_t)+2.$$
We check easily that $(\rho,\alpha)\in \{ \pm 3,\pm 6,\pm 9\}$, for $\alpha$ a long root. Hence the only possibility for $e(t)=\pm 1$ would come if
$(\de+\rho,\beta_t-\tilde\beta_t)=-3$. By observation (A), this forces
$(\rho,\beta_t-\tilde\beta_t)=-3$ and $(\lambda,\beta_t-\tilde\beta_t)=0$.
We deduce from this that
$\beta_t-\tilde\beta_t=(1,-2,1)$ and hence either $\beta_t=(1,-1,0)$,
$\tilde\beta_t=(0,1,-1)$, or $\tilde\beta_t=(-1,1,0)$,
$\beta_t=(0,-1,1)$. Moreover,
$$0=(\de,\beta_t-\tilde\beta_t)=\de_1-2\de_2+\de_3=-3\de_2.$$
But if $\delta_2=0$, $\delta+\beta_t$ would not be a dominant weight
in both cases.

{\it Case 2:}  $\beta_t\neq 0\neq \tilde\beta_t$ and $\beta_t+\tilde\beta_t$ is a long root.  In this case $(\beta_t,\tilde\beta_t)=1$,
$\beta_t-\tilde\beta_t$ is a short root and the formula for
$e(t)$ simplifies to
$$e(t)=(\de+\rho,\beta_t-\tilde\beta_t).$$
We check easily that $(\rho,\alpha)\in \{ \pm 1,\pm 4,\pm 5\}$, for $\alpha$ a short root. Again, using observation (A), it follows that
$(\rho,\beta_t-\tilde\beta_t)=\pm 1$, and $(\de,\beta_t-\tilde\beta_t)=0$.
It follows that $\beta_t-\tilde\beta_t=\pm (1,-1,0)$ and $\de_1=\de_2$.
But then there is only one path from $\delta$ to $\nu$, see e.g.
Remark \ref{tensorgeometry}

{\it Remaining cases:} If one of the weights $\beta_t$ or $\tilde\beta_t$ is equal to 0, the formula for $e(t)$ again simplifies to one of the two
versions in the two previous cases. The claim can be proven by similar
arguments as before. We leave this to the reader.

\begin{theorem}\label{surjectG2} Let $V$ be the 7-dimensional irreducible
representation of $G_2$, and let $V_\la$ be an irreducible representation
with highest weight $\la$. Then the image
of the affine braid group $AB_n$ generates $\End_\U(V_\la\otimes V^{\otimes n})$ for all $n\in\N$.
\end{theorem}

$Proof.$ The proof goes by induction on $n$. For $n=1$, we need to show
that $T$ generates $\End_\U(V_\la\otimes V)$. By Prop. \ref{prop:eigenvalue}, $T$ acts via the scalar $q^{2c_\mu-c_\la-c_V}$ on the submodule
$V_\mu\subset V_\la\otimes V$. The claim follows as soon as
we can show that $c_{\mu_1}\neq c_{\mu_2}$ for $\mu_1\neq\mu_2$.
Observe that if $V_{\mu_i}\subset V_\la\otimes V$, we can write $\mu_i=\la+\om_i$ for suitable weights $\om_i$ of $V$, for $i=1,2$.
Now
$$c_{\mu_1}-c_{\mu_2}= 2(\la+\rho,\om_1-\om_2)+(\om_1,\om_1)-(\om_2,\om_2).$$
If both $\om_1$ and $\om_2$ are nonzero, the second and third summands
on the right hand side cancel. The claim now follows from the fact
that $\om_1-\om_2$ is a root. One similarly shows the claim if
one of the $\om_i$'s is equal to 0, except possibly when $\la_1=\la_2$.
But in this case $V_\la\not\subset V_\la\otimes V$, so there is no second
path for which $T$ would have the same eigenvalue. This proves the
claim for $n=1$.

Observe that for the induction step, it suffices to show for each
module $W_\nu^{(n)}$ that the conditions of Lemma \ref{generatinglemma}
are satisfied. Using notations just in front of that lemma,
this means we have to show that $\mu_i$'s are equivalent
for any dominant weight $\mu_i$ for which $V_{\mu_i}\subset
V_\la\otimes V^{\otimes n-1}$ and for which $V_\nu\subset V_{\mu_i}\otimes V$.

In view of Prop. \ref{tensorold}, it suffices to show this for
the recent and new part of the tensor product.
Let $V_\nu$ be an irreducible representation in that part, and let
$V_\delta$ be any irreducible in $V_\la\otimes V^{\otimes n-2}$.
Then $\delta\neq\nu$, as otherwise $n(\nu)=n(\delta)\leq n-2$,
contradicting our assumption on $V_\nu$.
It remains  to deal with the following two cases,
as in the proof of Lemma \ref{checke(t)}. In order to simplify
notation, let us first assume that $\nu-\la$ is in the dominant Weyl
chamber.

{\it Case 1:}  $V_\nu$ is in the recent part, i.e. $n(\nu)=n-1$.
The idea is to show that the matrix entries of $P$ are non-zero
using Proposition \ref{AB2representations}.
We have
$$W^{(n)}_\nu\ \cong\ \bigoplus_\mu W^{(n-1)}_\mu,$$
where $\mu$ ranges over all dominant weights with $n(\mu)\leq n-1$
such that $V_\nu\subset V_\mu\otimes V$. It is easy to check,
using Example \ref{nnu},
that $\mu$ is contained in the set
$$\{ \nu,\nu\pm (\om_1-\om_2),\nu-\om_1,\nu-\om_2\},$$
where $\om_1$ and $\om_2$ are as in the proof of Lemma \ref{n(nu)}.
Geometrically, $\mu$ is either the center, or one of the
corners in the lower
half (including the separating line $\nu+t(\om_1-\om_2)$, $t\in\R$)
of the hexagon centered at $\nu$.
Taking $\delta=\nu-\om_1$, we get the equivalence
(in the sense used in Lemma \ref{generatinglemma}) of the weights
$\nu,\nu-\om_1,\nu-\om_2$ and $\nu-\om_1+\om_2$, by Prop. \ref{AB2representations} and Lemma \ref{checke(t)}.
The equivalence  of
$\nu+\om_1-\om_2$ with these weights is established in the same way,
for $\delta=\nu-\om_2$.
If $\nu$ is close to the boundary of the dominant Weyl chamber,
some of these $\mu$'s may not exist. However, the criteria for the
remaining paths would still be satisfied.

{\it Case 2:} If $n(\nu)=n$, we can only take the subset of $\mu$'s
from case 1 for which $n(\mu)=n-1$. These would be $\mu=\nu-\om_1$ and
$\mu=\nu-\om_2$, which forces $\delta=\nu-\om_1-\om_2$. Then we can
again use Lemma \ref{generatinglemma} together with Lemma \ref{AB2rep}
and Lemma \ref{checke(t)}.

{\it General case:} Let $W_\la$ be the reflection group generated by the hyperplanes going through $\la$ which are orthogonal to the roots.
Obviously, $W_\la\cong W$, the Weyl group. Moreover, let
$w\in W_\la$ be such that $w(\nu)$ is in the dominant Weyl chamber,
shifted by $\la$. Then the discussion in Cases 1 and 2 does apply
to $w(\nu)$.
Moreover, by Lemma \ref{n(nu)} we have $n(\nu)=n(w(\nu))$.
The combinatorics about the decomposition of $W^{(n)}_\nu$
can hence be deduced from the one of $W^{(n)}_{w(\nu)}$.
Essentially, the only difference is that we take the half of the
hexagon centered at $\nu$ which is closest to $\la$, in case 1.
All the other steps in the proof will be the same. Similar reasoning
also works in the much easier case 2.

\section{Structure Coefficients}

There exists an attractive graphical description of $End_\U(V^{\otimes n})$
in terms of spiders, see \cite{Ku2}. At this point, we do not have
a similarly satisfactory description of $\End_\U(V_\la\otimes V^{\otimes n})$. However, as a consequence of our approach, we are able to calculate
certain relations such as e.g. the one in Prop.  \ref{somerelations}, (c).

\subsection{Graphical Calculus}
We are going to use graphical calculus as in e.g. \cite{RT}, \cite{Turaev}
or \cite{Ks}. Here $R$-matrices are, as usual, given by crossing strands,
and the element $T$ is given by wrapping a string labeled by $V$ around
a string labeled by $V_\la$.
We also describe the mappings $i: V\to V^{\otimes 2}$ and $d:V^{\otimes 2}
\to V$ by the following trivalent graphs



\begin{center} \begin{tikzpicture}[scale=.6]


\path [ solid,black,thick,draw=none] (2.5, 10 cm) edge (2.5, 11 cm);
\path [ solid,black,thick,draw=none] (2.5, 11 cm) edge (3, 12 cm);
\path [ solid,black,thick,draw=none] (2.5, 11 cm) edge (2, 12 cm);

\path [ solid,black,thick,draw=none] (7.5, 10 cm) edge (8, 11 cm);
\path [ solid,black,thick,draw=none] (8.5, 10 cm) edge (8, 11 cm);
\path [ solid,black,thick,draw=none] (8, 11 cm) edge (8, 12 cm);

\node at (1,11 cm) {i=};
\node at (7,11 cm) {d=};

  \end{tikzpicture}
\end{center}
which are normalized such that $d\circ i=-([7]-1)\ 1_V= -(q^2+q^{-2})(q+1+q^{-1})1_V$.

\begin{proposition}\label{somerelations} We have the following relations:

(a) The eigenvalue of $T$ for $V_{\la+\om}\subset V_\la\otimes V$,
 with $\om$ a weight of $V$, is given by $q^{2(\la+\rho,\om)-10}$ for $\om\neq 0$,
 and by $q^{-12}$ for $\om=0$.

(b) The endomorphism of $V_\la$ given by the picture

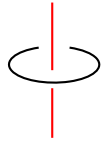
\begin{center}

 \begin{tikzpicture}[scale=.6]

\path [ solid,red,thick,draw=none] (2, 0 cm) edge (2, 1.5 cm);
\path [ solid,red,thick,draw=none] (2, -0.4 cm) edge (2, -1.5 cm);

\draw[thick, black] (1.7,0.5) arc (110:430:1 and 0.4);

  \end{tikzpicture}

  \end{center}
is given by the following scalar, where the summation goes over the elements of the Weyl group
$$d=\frac{\sum_w \varepsilon(w) q^{(\la+\rho, w(\La_1+\rho))}}
{\sum_w \varepsilon(w) q^{(\la+\rho, w(\rho))}}.$$

(c) The eigenvalue of the endomorphism $(1\otimes d)(R_{V,V_\la}R_{V_\la,V}\otimes 1)(1\otimes i)$ of $V_\la\otimes V$
for $V_\nu\subset V_\la\otimes V$, $\nu\neq\la$,  given by the picture below is equal to
 $c(\la,\nu)$, as in Lemma \ref{clambdanu}.

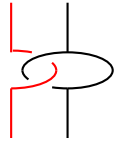
\begin{center}

 \begin{tikzpicture}[scale=.6]

\path [ solid,black,thick,draw=none] (2, 0.4 cm) edge (2, 1.5 cm);
\path [ solid,black,thick,draw=none] (2, -0.4 cm) edge (2, -1.5 cm);
\path [ solid,red,thick,draw=none] (0.75, 0.4 cm) edge (0.75, 1.5 cm);
\path [ solid,red,thick,draw=none] (0.75, -0.4 cm) edge (0.75, -1.5 cm);

\draw[thick] (2,-0.4) arc (-90:180:1 and 0.4);
\draw[thick] (1,0) arc (-180:-150:1 and 0.4);
\draw[thick] (2,-0.4) arc (-90:-110:1 and 0.4);

\draw[thick, red] (1.75,0) arc (0:-90:1 and 0.4);
\draw[thick, red] (1.75,0) arc (0:25:1 and 0.4);
\draw[thick, red] (1.21,0.4) arc (65:90:1 and 0.4);

  \end{tikzpicture}
\end{center}
\end{proposition}

\medskip

$Proof.$ Claim (a) follows from Prop. \ref{prop:eigenvalue}. Statement (b) is well-known in quantum topology. See e.g. \cite{Turaev}, \cite{RT} and for an explicit calculation of the entries of the $S$-matrix e.g. \cite{TW}.
Statement (c) will follow from Lemma \ref{clambdanu}, which will be proved
in the following sections.

\subsection{Matrix coefficients of $P$} We consider representations of the affine braid group $AB_2$ as in Section \ref{AB2section}. In particular, we denote by $P$ the  rank 1 eigenprojection belonging to the third eigenvalue $r^{-1}=-q^{-7}$ of $A$, and $T$ is given as a diagonal matrix
$T={\rm diag}(q^{e(t)})$, with $e(t)$ as in the proof of Lemma
\ref{checke(t)}; this differs from the values in Prop. \ref{somerelations}
by a constant factor for each representation $W=W_\nu=\Hom(V_\nu,V_\lambda\otimes V^{\otimes 2})$.
However, this does not have any effect on the result we need.
The following lemma can be considered a version
of \cite{Wexc}, Prop 5.7.

\begin{lemma}\label{diagonalentries} Let $d_s$ be the diagonal entry of $P$
with respect to the path $s$. Assume that $\prod_t q^{e(t)}=-rq^{-1}$.
Then $d_s$ is given by
$$d_s=\frac{[e(s)+1]}{1- [7]}\ \prod_{t\neq s}
\frac{[(e(s)+e(t))/2]}{ [(e(s)-e(t))/2]}.$$
\end{lemma}

$Proof.$ One proceeds as in the proof of \cite{Wexc}, Prop 5.7.
Basically, we express the coefficients of $A$ in terms of the ones
of $P$, see Prop. \ref{AB2representations}, and plug this into
the equation $AP=r^{-1}P$, where $r=-q^7$, see Example
\ref{G2parameters}.  We only have
to make the following adjustment: It was claimed there that
$\prod q^{e(t)}=rq$.
However, looking at the proof more carefully, one sees that there is a second possibility, namely $\prod_t q^{e(t)}=-rq^{-1}$ (indeed, the equation
$d_s-\bar d_s=0$ gives a quadratic equation for $\prod_t q^{e(t)}$
with two solutions). Proceeding under this
assumption as in that proof, we obtain the formula for $d_s$ as claimed in the statement.

\subsection{$G_2$-representations} We now apply the results from the previous section to the representations in connection with $G_2$.
Let $\beta,\beta_1$ and $\beta_2$
be nonzero weights of the 7-dimensional representation $V$ such that
$\beta_1+\beta_2=\beta$. If we set $\nu=\lambda+\beta$, we have four paths of length 2 from $\beta$ to $\nu$, with the middle weight $\mu_t$
running through $\lambda,\lambda+\beta_1,\lambda+\beta_2$ and $\lambda+\beta$.
 We define the scalar $c(\lambda,\beta)$ by
$(1-[7])\ PTP=c(\lambda,\beta)P$. As $T={\rm diag}(q^{f(t)})$,
with $f(t)$ as given in Prop. \ref{somerelations},(a),
$$c(\lambda,\nu)=(1-[7])\ \sum_t q^{f(t)}d_t.$$

\begin{lemma}\label{clambdanu}
The coefficient $c(\la,\nu)$ is given by
$$c(\la,\nu)=q^{-5}\ [(q+q^{-1})(q^{2(\la+\rho,\beta_1)}+q^{2(\la+\rho,\beta_2)})
+ q^{2(\la+\rho,\beta)+1}+q^{-1}].$$
\end{lemma}

$Proof.$
Using the second formula for $e(t)$ in the proof of Lemma \ref{checke(t)},
one checks that
$\sum e(t)=6$; indeed, it suffices to observe that interchanging
$\beta_t$ with $\tilde\beta_t$ describes another path in that sum.
Using this and $r=-q^7$,  it follows from Lemma \ref{diagonalentries} that
$$(1-[7])\ \sum_s q^{e(s)}d_s\ =\ \sum_s q^{e(s)}
[e(s)+1]\ \prod_{t\neq s}
\frac{[(e(s)+e(t))/2]}{ [(e(s)-e(t))/2]}.
$$
 We obtain, aided by computer software, that this is equal to

$$q^7\  \sum_t [e(t)].$$
Observe that by Prop \ref{prop:eigenvalue} or Prop. \ref{somerelations},(a),
the eigenvalues $q^{f(t)}$
of $T$ are given by $f(t)=c_{\mu_t}-c_\la-c_V$. Hence, in order to get
$c(\la,\nu)$, we have to multiply the expression above by
$q^{f(t)-e(t)}= q^{(c_\nu-c_\la)/2-c_V-1}=q^{(\la+\rho,\beta)-12}$.
We obtain the claimed result after simplifying the sum,
using the second formula for $e(t)$ in the proof of Lemma \ref{checke(t)}.

\begin{remark} We can express all the constants in the relations in
this section as rational functions in $q$, $r_1$ and $r_2$,
where $r_i=q^{\la_i}$, $i=1,2$. This can be seen either by
interpolating the formulas for integral dominant weights $\la$,
or doing the same calculations for arbitrary weights $\la$ with $V_\la$ being the Verma module for that weight. It would be desirable to define
a generic quotient of the affine braid group just in terms of relations
involving these variables, independently of $\la$.
\end{remark}
\subsection{General discussion}
Our method was based on the fact that braid representations for which
the generators only have three distinct eigenvalues can be studied
via an extension of the Jucys-Murphy approach. So under suitable conditions
the tensor product rules essentially determine the braid representation.
We give here a brief informal discussion how this can be used for other
exceptional Lie types:

(a) The $E_N$ series in \cite{Wexc}. Here the module $V$ was chosen to be
the fundamental module corresponding to the vertex in the Dynkin diagram
$E_N$  furthest from the triple point for $N>5$, $N\neq 9$.
Path representations of braid groups were defined for a combinatorially
defined summand $V^{\otimes n}_{new}\subset V^{\otimes n}$. For $N<9$
this indeed corresponds to the part of the tensor product as defined
in this paper. One can deduce from this a two-parameter family of
braid representations which contain the Hecke algebra and $BMW$-algebra
representations as quotients. As an application, it was shown that
for $N=6,7$ the image of the braid group $B_n$ and one additional
element generate $\End_\U(V^{\otimes n})$, with $\U=U_q\g(E_N)$.

(b) We have made additional studies regarding representations of
$F_4$. If $V$ is the representation of dimension 26, one can show
that the braid representation into $\End_\U(V^{\otimes n}_{new})$ factors
through the $BMW$-algebra quotient for the vector representation of
the symplectic quantum group $U_q\g(C_3)$
(see \cite{BBMW}, \cite{Mu} and \cite{WBMW}).
If one takes for $V$ the
adjoint representation, one again obtains braid representations
into $\End_\U(V^{\otimes n}_{new})$ for which the generators only
have three eigenvalues. They seem to be related to the $E_N$-series,
but more complicated.

(c) One of the motivations for this study was the proposed exceptional
series of Vogel and Deligne. These authors observed a certain
uniform behavior for the decomposition of $small$
tensor powers of the adjoint
representations of Lie algebras of exceptional Lie type.
Here again the braid representations in the new part of
these tensor powers have the property that the
braid generators satisfy a cubic equation. We expect that
techniques in \cite{Wexc} and in this paper should be useful
to establish a uniform behavior for a significant part
of arbitrarily large tensor powers.

\end{document}